\documentclass[11pt,twoside,a4paper]{article}

\usepackage{verbatim}
\usepackage{amsmath}
\usepackage{amsfonts}
\usepackage{url}
 
\usepackage{pgf,tikz}

\newcommand\polygon[3][]{
  \pgfmathsetmacro{\angle}{360/#2}
  \pgfmathsetmacro{\polyangle}{180*(#2-2)/#2}  
  \pgfmathsetmacro{\startangle}{-180+\polyangle/2}
  \pgfmathsetmacro{\y}{cos(\angle/2)}
  \begin{scope}[#1]

    \node at (0,0) {\textbullet};
  
    \foreach \i in {1,2,...,#2} 
    {
      \pgfmathsetmacro{\x}{\startangle - \angle*\i}
      
      \draw[fill] (\x:\y cm) -- (\x - \angle:\y cm);    
      
    }

	\foreach [var=\source,var=\target] in #3 
	{
	    \pgfmathsetmacro{\x}{\startangle+ \angle}	

		\draw (\x  - \angle*\source : \y cm) -- (\x - \angle*\target : \y cm);		      
	
	}    
	    
  \end{scope}
}

\newcommand\lpolygon[3][]{
  \pgfmathsetmacro{\angle}{360/#2}
  \pgfmathsetmacro{\polyangle}{180*(#2-2)/#2}  
  \pgfmathsetmacro{\startangle}{-180+\polyangle/2}
  \pgfmathsetmacro{\y}{cos(\angle/2)}
  \begin{scope}[#1]

    \node at (0,0) {\textbullet};
  
    \foreach \i in {1,2,...,#2} 
    {
      \pgfmathsetmacro{\x}{\startangle - \angle*\i}
      
      \draw[fill] (\x:\y cm) -- (\x - \angle:\y cm);    
      
      \draw (\x+ \angle+ 4:4+\y cm) node{\i};
      
    }

	\foreach [var=\source,var=\target] in #3 
	{
	    \pgfmathsetmacro{\x}{\startangle+ \angle}	

		\draw (\x  - \angle*\source : \y cm) -- (\x - \angle*\target : \y cm);		      
	
	}    
	    
  \end{scope}
}

\definecolor{zzttgg}{rgb}{0.0,0.8,0}

  \title{Again around frieze patterns}
  \author{Tiberiu Spircu \\ Institute of Mathematics of the Romanian Academy \\ 21, Calea Grivi\c{t}ei. 010736 Bucharest, Romania \\ \it email: spircut@yahoo.com \rm \\
\small and \\  
 Stefan V. Pantazi\\
Conestoga College ITAL\\
299 Doon Valley Drive, Kitchener ON, Canada N2G 4M4\
\\ \it email: svpantazi@gmail.com \rm }
 
   \date{\today}
\begin{document}

\maketitle
 
\begin{abstract}
The main goal of this paper is to prove several new results about frieze patterns and their equivalents, the quiddity (or $\eta$-)sequences and to obtain a formula giving the number of non-similar frieze patterns of given finite width.
\end{abstract}  

 \it AMS classification\rm: 05A17, 05E18.

\vspace{3pt}
 \it Keywords\rm: frieze pattern, quiddity sequence, $\eta$-sequence, special linear group.

  \section{Introduction}

Frieze patterns have a long history. Starting from the ``immortals archers" in Persepolis, continuing with wall decorations in Alhambra, more recently with intriguing drawings of Maurits Cornelis Escher, they have raised many questions to specific-shaped minds.

Since 1977, when the paper of Conway and Coxeter (\cite{ConwayCoxeter}) appeared, in which a link between particular finite frieze patterns of numbers and triangulations of convex polygons was first established, many other connections with other mathematical objects have been found (see \cite{ConwayCoxeter}).

In this paper two kinds of new results about frieze patterns are presented: first (Section 2) a link with elements of finite index in the special linear group $SL_2(\mathbb{Z})$; second, a method to calculate non-similar frieze patterns (Section 4). The main result in Section 3 concerns the so-called basic sequences: they can be supplemented to quiddity sequences.

  \section{Definitions and notations}

The special linear group $SL_2(\mathbb{Z})$ is well-known \cite{Conrad,Humphreys}. Its elements are integer $2\times 2$ matrices $A=\left( \begin{array}{cc} a & b \\ c & d \end{array} \right)$ whose determinant $a d - b c$ is equal to 1. This group is generated (see for example \cite{Humphreys}, p.80) by two matrices, $S=\left( \begin{array}{cc} 0 & 1 \\ -1 & 0 \end{array} \right)$ and $T=\left( \begin{array}{cc} 0 & 1 \\ -1 & -1 \end{array} \right)$, of order 4, res. 3; in fact, it is the direct product of the multiplicative cyclic group $\{-, +\}\simeq \mathbb{Z}_2$ with the free product of $ \mathbb{Z}_2$ and $ \mathbb{Z}_3$ (see \cite{Alperin}). In \cite{Humphreys} it is shown that any element $A\in SL_2(\mathbb{Z})$ is expressed as follows

\vspace{3 pt}
$\pm T^{b_0}ST^{e_1}S...T^{e_n}S^{b_1}$

\vspace{3 pt}
\noindent (where $b_0, b_1\in\{0, 1\}$ and $e_i\in \{1, 2\}$) as a product of generators and any such expression never reduces to $\pm I$.

Unfortunately, the matrix $T$ has no ``geometric" interpretation. Much more useful is the ``translation" matrix $U=\left( \begin{array}{cc} 1 & 0 \\ 1 & 1 \end{array} \right)$, which is expressed as $ST^2$. Since $U^a=\left( \begin{array}{cc} 1 & 0 \\ a & 1 \end{array} \right)$, the matrix $U$ has infinite order.

It is easy to check that $US$ has order 6, but $U^a S$, for $a>1$, has infinite order; moreover, $USU^2S$ has order 4 and $U^{-1}S=SUSU$ has order 3.

What can be said about the expression $S^{b_0}U^{a_1}SU^{a_2}S\cdot\cdot\cdot U^{a_n}S^{b_1}$ (where $b_0, b_1\in\{0, 1\}$ and $a_i\ge 1$)? Is its order finite or infinite? Does it reduce to $\pm I$ or not?

It is enough to study the expressions $U^{a_1}SU^{a_2}S\cdot\cdot\cdot U^{a_n}S$, since if such an expression $A$ reduces to $\pm I$, then both $SA$ and $AS$ reduce to $\pm S$, thus $(AS)^2=(SA)^2=-I$, and $SAS$ reduces to $\mp I$.


Note that if the expression $U^{a_1}SU^{a_2}S\cdot\cdot\cdot U^{a_n}S$ reduces to $\pm I$, then the ``rotated" expression  $U^{a_2}SU^{a_3}S...U^{a_n}SU^{a_1}S$ also reduces to $\pm I$.

In addition, $(U^{a+1}S)(US)(U^{b+1}S)=(U^{a}S)(U^{b}S)$ for $a, b \in \mathbb{Z}$. The relations above suggest considering $\eta$-sequences.

\vspace{3pt}
\bf Definition 1 \rm (see \cite{CuntzHeckenberger}, Definition 3.2). $\eta$\bf-sequences \rm (of positive integers) of length $n\ge 3$ are defined recursively as follows:

\vspace{3pt}
\noindent (Rule 1 -- initialize) $(1,1,1)$ is an $\eta$-sequence of length $1$.

\vspace{3pt}
\noindent (Rule 2 -- rotate) If $(c_0,c_1,...,c_{n-2},c_{n-1})$ is an $\eta$-sequence of length $n$, then 
\indent $(c_1,c_2,...,c_{n-1},c_0)$ is an $\eta$-sequence of same length.

\vspace{3pt}
\noindent (Rule 3 -- expand) If $(c_0,c_1,c_2,...,c_{n-1})$ is an $\eta$-sequence of length $n$, then  
\indent $(c_0+1,1,c_1+1,c_2,...,c_{n-1})$ is an $\eta$-sequence of length $n+1$.

\vspace{6pt}
An easy reasoning shows that

\vspace{3pt}
\noindent (Rule 2' -- reverse) If $(c_0,c_1,...,c_{n-2},c_{n-1})$ is an $\eta$-sequence of length $n$, then  

\indent $(c_{n-1},c_{n-2},...,c_1,c_0)$ is an $\eta$-sequence of same length

\vspace{3pt}
\noindent is a direct consequence of the rules 1, 2, 3 above.

Thus, starting from the relation $(US)^3=-I$ which corresponds to the initial $\eta$-sequence $(1,1,1)$, the following is immediate:

\bf Proposition 1\rm . For \ any \ $\eta$-sequence \  $(c_0,c_1,...,c_{n-1})$, \ the \ expression $U^{c_0}SU^{c_1}S...U^{c_{n-1}}S$ reduces to $-I$.


\vspace{6pt}
Let us end this section by presenting the following easy-to-prove:

\bf Lemma 1\rm . If $X\in SL_2(\mathbb{Z})$ and $a, b$ are two integers such that $XU^{a}S =U^{b}SX$, then $a=b$. $\square$

\vspace{6pt}

\subsection {Dual representation}

Is is well known (see Remark 2.4 from \cite{Cuntz}) that any $\eta$-sequence of length $n$ corresponds to a triangulation of a convex $n$-gon by non-intersecting diagonals; namely, the component $c_i$ of the sequence is exactly the number of triangles supported by vertex $i$.

It is also known (\cite{O'Rourke}, Lemma 1.3, pp. 12-13) that the ``weak" dual graph of an $n$-gon triangulation is a tree with each node of degree at most 3. For our purposes, the dual graph is a \it full \rm binary tree with $n-2$ internal nodes corresponding to triangles, a fixed root labeled $\bf a$ and $n-1$ leaves labeled $b,c,d...$ associated with the other sides of the $n$-gon. The count of consecutive internal (i.e., non-leaf) nodes visited by a depth-first traversal of the binary tree yields the component $c_i$ of an $\eta$-sequence. For example, the pentagon diagonalization corresponding to the $\eta$-sequence 12213 and its dual are illustrated in the following. 

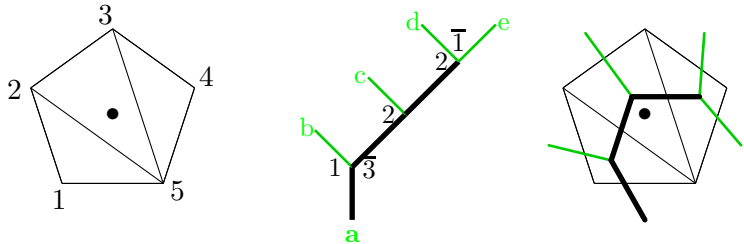
\begin{tikzpicture}[scale=.7]

\lpolygon[scale=2,xshift=0cm,yshift=5cm]{5}{{3/5,2/5}}

\begin{scope} [scale=1,xshift=2cm,yshift=3cm]
\draw (2.5,5)		node[anchor=north] {\small\color{green}\bf a};
\draw [line width=2pt] (2.5,5)-- (2.5,6);
\draw [color=zzttgg, line width=1pt] (2.5,6)-- (1.8,6.7);
\draw (1.65,6.4)		node[anchor=south] {\small\color{green}b};
\draw (2.5,6)		node[anchor=east] {\small1};
\draw [line width=2pt] (2.5,6)-- (3.5,7);
\draw [color=zzttgg, line width=1pt] (3.5,7)-- (2.8,7.7);
\draw (2.65,7.4)		node[anchor=south] {\small\color{green}c};
\draw (3.5,7)		node[anchor=east] {\small2};
\draw [line width=2pt] (3.5,7)-- (4.5,8);
\draw [color=zzttgg, line width=1pt] (4.5,8)-- (3.8,8.7);
\draw (3.65,8.4)		node[anchor=south] {\small\color{green}d};
\draw (4.5,8)		node[anchor=east] {\small2};
\draw [color=zzttgg, line width=1pt] (4.5,8)-- (5.2,8.7);
\draw (5.35,8.4)		node[anchor=south] {\small\color{green}e};
\draw (4.5,8)		node[anchor=south] {\small$\overline{1}$};
\draw (2.5,6)		node[anchor=west] {\small$\overline{3}$};
\end{scope}

\begin{scope}[scale=1,xshift=0cm,yshift=0cm]

\polygon[scale=2,xshift=5cm,yshift=5cm]{5}{{3/5,2/5}}

\draw [color=zzttgg, line cap=round, line join=round, line width=1pt] (11.026,10.333)-- (11.117,11.537);
\draw [color=zzttgg, line cap=round, line join=round, line width=1pt] (11.026,10.333)-- (11.807,9.413);
\draw [color=zzttgg, line cap=round, line join=round, line width=1pt] (9.758,10.333)-- (8.883,11.537);
\draw [line cap=round, line join=round, line width=2pt] (9.758,10.333)-- (11.026,10.333);
\draw [color=zzttgg, line cap=round, line join=round, line width=1pt] (9.366,9.127)-- (8.193,9.413);
\draw [line cap=round, line join=round, line width=2pt] (9.366,9.127)-- (9.758,10.333);
\draw [line cap=round, line join=round, line width=2pt] (9.366,9.127)-- (10.000,8.000);

\end{scope}

\end{tikzpicture}


The choices for the graphical representation of the dual graph underline the observation in \cite{Morier-Genoud} that Coxeter's frieze patterns are indeed at the crossroads of algebra, geometry and combinatorics. A simple, straight angles approach often suffices for illustrating  algebraic or combinatorial properties of the binary trees. However, a triangulated $n$-gon overlay (on the right) that prevents the overlap of branches for any arbitrary tree, may render geometric representations more ``natural".

From the dual graph point of view, rules 1-3 defining $\eta$-sequences are  interpreted as follows. Obviously, Rule 1 defines the smallest non-empty, \it full \rm binary tree (1,1,1).

Rule 2 defines the rotation and reverse rotation of a binary tree such as (1,3,2,1,5,1,2,3) $\leftrightarrow$ (3,2,1,5,1,2,3,1)

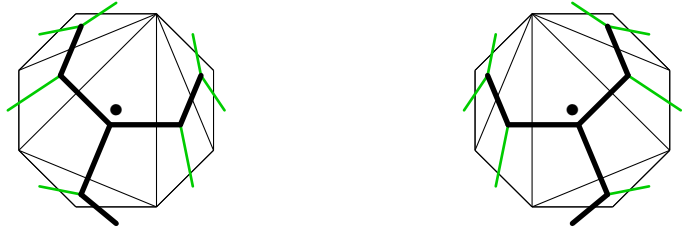
\begin{tikzpicture}[scale=.75]

\begin{scope}[scale=1,xshift=0cm,yshift=0cm]

\polygon[scale=2,xshift=5cm,yshift=5cm]{8}{{3/5,2/5,5/7,5/8,2/8}}
\draw [color=zzttgg, line cap=round, line join=round, line width=1pt] (9.384,11.487)-- (8.656,11.344);
\draw [color=zzttgg, line cap=round, line join=round, line width=1pt] (9.384,11.487)-- (10.000,11.900);
\draw [color=zzttgg, line cap=round, line join=round, line width=1pt] (9.023,10.616)-- (8.100,10.000);
\draw [line cap=round, line join=round, line width=2pt] (9.023,10.616)-- (9.384,11.487);
\draw [color=zzttgg, line cap=round, line join=round, line width=1pt] (11.487,10.616)-- (11.344,11.344);
\draw [color=zzttgg, line cap=round, line join=round, line width=1pt] (11.487,10.616)-- (11.900,10.000);
\draw [line cap=round, line join=round, line width=2pt] (11.126,9.745)-- (11.487,10.616);
\draw [color=zzttgg, line cap=round, line join=round, line width=1pt] (11.126,9.745)-- (11.344,8.656);
\draw [line cap=round, line join=round, line width=2pt] (9.894,9.745)-- (9.023,10.616);
\draw [line cap=round, line join=round, line width=2pt] (9.894,9.745)-- (11.126,9.745);
\draw [color=zzttgg, line cap=round, line join=round, line width=1pt] (9.384,8.513)-- (8.656,8.656);
\draw [line cap=round, line join=round, line width=2pt] (9.384,8.513)-- (9.894,9.745);
\draw [line cap=round, line join=round, line width=2pt] (9.384,8.513)-- (10.000,8.000);

\end{scope}

\begin{scope}[scale=1,xshift=8cm,yshift=0cm]

\polygon[scale=2,xshift=5cm,yshift=5cm]{8}{{2/4,1/4,4/6,4/7,1/7}}
\draw [color=zzttgg, line cap=round, line join=round, line width=1pt] (8.513,10.616)-- (8.100,10.000);
\draw [color=zzttgg, line cap=round, line join=round, line width=1pt] (8.513,10.616)-- (8.656,11.344);
\draw [color=zzttgg, line cap=round, line join=round, line width=1pt] (8.874,9.745)-- (8.656,8.656);
\draw [line cap=round, line join=round, line width=2pt] (8.874,9.745)-- (8.513,10.616);
\draw [color=zzttgg, line cap=round, line join=round, line width=1pt] (10.616,11.487)-- (10.000,11.900);
\draw [color=zzttgg, line cap=round, line join=round, line width=1pt] (10.616,11.487)-- (11.344,11.344);
\draw [line cap=round, line join=round, line width=2pt] (10.977,10.616)-- (10.616,11.487);
\draw [color=zzttgg, line cap=round, line join=round, line width=1pt] (10.977,10.616)-- (11.900,10.000);
\draw [line cap=round, line join=round, line width=2pt] (10.106,9.745)-- (8.874,9.745);
\draw [line cap=round, line join=round, line width=2pt] (10.106,9.745)-- (10.977,10.616);
\draw [line cap=round, line join=round, line width=2pt] (10.616,8.513)-- (10.106,9.745);
\draw [color=zzttgg, line cap=round, line join=round, line width=1pt] (10.616,8.513)-- (11.344,8.656);
\draw [line cap=round, line join=round, line width=2pt] (10.616,8.513)-- (10.000,8.000);

\end{scope}

\end{tikzpicture}\\


Rule 3 allows for the left sided expansion of a binary tree with $n$ leaves  by growing a new first internal node (just above root) and a right leaf node. The following illustrates the expansion (1,1,1) $\rightarrow$ (2,1,2,1) $\rightarrow$ (3,1,2,2,1) $\rightarrow$ (4,1,2,2,2,1).

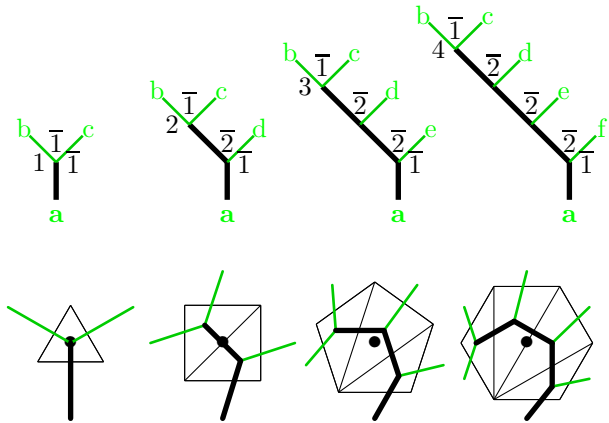
\begin{tikzpicture}[scale=.5]

\begin{scope}[scale=1,xshift=0cm,yshift=3cm]

\draw (1.5,3)		node[anchor=north] {\small\color{green}\bf a};
\draw [line width=2pt] (1.5,3)-- (1.5,4);
\draw [color=zzttgg, line width=1pt] (1.5,4)-- (0.8,4.7);
\draw (0.65,4.4)		node[anchor=south] {\small\color{green}b};
\draw (1.5,4)		node[anchor=east] {\small1};
\draw [color=zzttgg, line width=1pt] (1.5,4)-- (2.2,4.7);
\draw (2.35,4.4)		node[anchor=south] {\small\color{green}c};
\draw (1.5,4)		node[anchor=south] {\small$\overline{1}$};
\draw (1.5,4)		node[anchor=west] {\small$\overline{1}$};

\end{scope}

\begin{scope}[scale=1,xshift=4cm,yshift=2cm]

\draw (2,4)		node[anchor=north] {\small\color{green}\bf a};
\draw [line width=2pt] (2,4)-- (2,5);
\draw [line width=2pt] (2,5)-- (1,6);
\draw [color=zzttgg, line width=1pt] (1,6)-- (0.3,6.7);
\draw (0.15,6.4)		node[anchor=south] {\small\color{green}b};
\draw (1,6)		node[anchor=east] {\small2};
\draw [color=zzttgg, line width=1pt] (1,6)-- (1.7,6.7);
\draw (1.85,6.4)		node[anchor=south] {\small\color{green}c};
\draw (1,6)		node[anchor=south] {\small$\overline{1}$};
\draw [color=zzttgg, line width=1pt] (2,5)-- (2.7,5.7);
\draw (2.85,5.4)		node[anchor=south] {\small\color{green}d};
\draw (2,5)		node[anchor=south] {\small$\overline{2}$};
\draw (2,5)		node[anchor=west] {\small$\overline{1}$};

\end{scope}

\begin{scope}[scale=1,xshift=8cm,yshift=1cm]

\draw (2.5,5)		node[anchor=north] {\small\color{green}\bf a};
\draw [line width=2pt] (2.5,5)-- (2.5,6);
\draw [line width=2pt] (2.5,6)-- (1.5,7);
\draw [line width=2pt] (1.5,7)-- (0.5,8);
\draw [color=zzttgg, line width=1pt] (0.5,8)-- (-0.2,8.7);
\draw (-0.35,8.4)		node[anchor=south] {\small\color{green}b};
\draw (0.5,8)		node[anchor=east] {\small3};
\draw [color=zzttgg, line width=1pt] (0.5,8)-- (1.2,8.7);
\draw (1.35,8.4)		node[anchor=south] {\small\color{green}c};
\draw (0.5,8)		node[anchor=south] {\small$\overline{1}$};
\draw [color=zzttgg, line width=1pt] (1.5,7)-- (2.2,7.7);
\draw (2.35,7.4)		node[anchor=south] {\small\color{green}d};
\draw (1.5,7)		node[anchor=south] {\small$\overline{2}$};
\draw [color=zzttgg, line width=1pt] (2.5,6)-- (3.2,6.7);
\draw (3.35,6.4)		node[anchor=south] {\small\color{green}e};
\draw (2.5,6)		node[anchor=south] {\small$\overline{2}$};
\draw (2.5,6)		node[anchor=west] {\small$\overline{1}$};

\end{scope}

\begin{scope}[scale=1,xshift=12cm,yshift=0cm]

\draw (3,6)		node[anchor=north] {\small\color{green}\bf a};
\draw [line width=2pt] (3,6)-- (3,7);
\draw [line width=2pt] (3,7)-- (2,8);
\draw [line width=2pt] (2,8)-- (1,9);
\draw [line width=2pt] (1,9)-- (0,10);
\draw [color=zzttgg, line width=1pt] (0,10)-- (-0.7,10.7);
\draw (-0.85,10.4)		node[anchor=south] {\small\color{green}b};
\draw (0,10)		node[anchor=east] {\small4};
\draw [color=zzttgg, line width=1pt] (0,10)-- (0.7,10.7);
\draw (0.85,10.4)		node[anchor=south] {\small\color{green}c};
\draw (0,10)		node[anchor=south] {\small$\overline{1}$};
\draw [color=zzttgg, line width=1pt] (1,9)-- (1.7,9.7);
\draw (1.85,9.4)		node[anchor=south] {\small\color{green}d};
\draw (1,9)		node[anchor=south] {\small$\overline{2}$};
\draw [color=zzttgg, line width=1pt] (2,8)-- (2.7,8.7);
\draw (2.85,8.4)		node[anchor=south] {\small\color{green}e};
\draw (2,8)		node[anchor=south] {\small$\overline{2}$};
\draw [color=zzttgg, line width=1pt] (3,7)-- (3.7,7.7);
\draw (3.85,7.4)		node[anchor=south] {\small\color{green}f};
\draw (3,7)		node[anchor=south] {\small$\overline{2}$};
\draw (3,7)		node[anchor=west] {\small$\overline{1}$};

\end{scope}

\end{tikzpicture}\\

\begin{tikzpicture}[scale=.5]

\begin{scope}[scale=1,xshift=0cm,yshift=0cm]

\polygon[scale=2,xshift=5cm,yshift=5cm]{3}{{}}
\draw [color=zzttgg, line cap=round, line join=round, line width=1pt] (10.000,10.000)-- (8.355,10.950);
\draw [color=zzttgg, line cap=round, line join=round, line width=1pt] (10.000,10.000)-- (11.645,10.950);
\draw [line cap=round, line join=round, line width=2pt] (10.000,10.000)-- (10.000,8.000);

\end{scope}

\begin{scope}[scale=1,xshift=4cm,yshift=0cm]

\polygon[scale=2,xshift=5cm,yshift=5cm]{4}{{1/3}}
\draw [color=zzttgg, line cap=round, line join=round, line width=1pt] (9.529,10.471)-- (8.100,10.000);
\draw [color=zzttgg, line cap=round, line join=round, line width=1pt] (9.529,10.471)-- (10.000,11.900);
\draw [line cap=round, line join=round, line width=2pt] (10.471,9.529)-- (9.529,10.471);
\draw [color=zzttgg, line cap=round, line join=round, line width=1pt] (10.471,9.529)-- (11.900,10.000);
\draw [line cap=round, line join=round, line width=2pt] (10.471,9.529)-- (10.000,8.000);

\end{scope}

\begin{scope}[scale=1,xshift=8cm,yshift=0cm]

\polygon[scale=2,xshift=5cm,yshift=5cm]{5}{{1/3,1/4}}
\draw [color=zzttgg, line cap=round, line join=round, line width=1pt] (8.974,10.333)-- (8.193,9.413);
\draw [color=zzttgg, line cap=round, line join=round, line width=1pt] (8.974,10.333)-- (8.883,11.537);
\draw [line cap=round, line join=round, line width=2pt] (10.242,10.333)-- (8.974,10.333);
\draw [color=zzttgg, line cap=round, line join=round, line width=1pt] (10.242,10.333)-- (11.117,11.537);
\draw [line cap=round, line join=round, line width=2pt] (10.634,9.127)-- (10.242,10.333);
\draw [color=zzttgg, line cap=round, line join=round, line width=1pt] (10.634,9.127)-- (11.807,9.413);
\draw [line cap=round, line join=round, line width=2pt] (10.634,9.127)-- (10.000,8.000);

\end{scope}

\begin{scope}[scale=1,xshift=12cm,yshift=0cm]

\polygon[scale=2,xshift=5cm,yshift=5cm]{6}{{1/3,1/4,1/5}}
\draw [color=zzttgg, line cap=round, line join=round, line width=1pt] (8.667,10.000)-- (8.355,9.050);
\draw [color=zzttgg, line cap=round, line join=round, line width=1pt] (8.667,10.000)-- (8.355,10.950);
\draw [line cap=round, line join=round, line width=2pt] (9.667,10.577)-- (8.667,10.000);
\draw [color=zzttgg, line cap=round, line join=round, line width=1pt] (9.667,10.577)-- (10.000,11.900);
\draw [line cap=round, line join=round, line width=2pt] (10.667,10.000)-- (9.667,10.577);
\draw [color=zzttgg, line cap=round, line join=round, line width=1pt] (10.667,10.000)-- (11.645,10.950);
\draw [line cap=round, line join=round, line width=2pt] (10.667,8.845)-- (10.667,10.000);
\draw [color=zzttgg, line cap=round, line join=round, line width=1pt] (10.667,8.845)-- (11.645,9.050);
\draw [line cap=round, line join=round, line width=2pt] (10.667,8.845)-- (10.000,8.000);

\end{scope}

\end{tikzpicture}\\

\subsection {Frieze patterns and quiddity sequences}

\vspace{6 pt}
\bf Definition 2  \rm (see \cite{Coxeter}). A \bf frieze pattern (of integers) \rm is defined as a double-indexed sequence of integers
$\{\varphi (i,j)\}_{i\in  \mathbb{N}, j\in  \mathbb{Z}}$
satisfying the following two conditions:

\vspace{3pt}
\noindent (FI1) (initial conditions) $\varphi (0,j)=0$ and $\varphi (1,j)=1$ for each $j\in  \mathbb{Z}$; 

\vspace{3pt}
\noindent (FI2) (diamond rule) $\varphi (i,j)=\frac{\varphi (i-1,j+1) \varphi (i-1,j)-1}{\varphi (i-2,j+1)}$  for each $i\ge 3, j\in  \mathbb{Z}$.

\vspace{6pt}
Note that in a frieze pattern of integers all values in rows 3, 4, ... are well determined provided the 2-nd row is known.

\vspace{6 pt}
\bf Definition 3\rm . According to Conway-Coxeter (see \cite{ConwayCoxeter}), a $\bf quiddity\ sequence \rm$  is a finite sequence of positive integers

\vspace{3pt}
$(a_0,a_1,\dots,a_{n-1}) \ \ \ \ (n\ge 3)$

\vspace{3pt}
\noindent such that for the (non-trivial) frieze pattern in which 

\vspace{3 pt}
$\varphi (2,j)=a_{j\%n}  \rm \ for\ all\ \it  j\in \mathbb{Z}$

\vspace{3 pt}
\noindent there exists a ``row" $r\ge 3$ such that $\varphi (r,j)=1\ $ for all $j\in \mathbb{Z}$. (Of course, \% denotes here the remainder in an integer division.)

\vspace{6 pt}
In \cite{ConwayCoxeter} it was also established that any quiddity sequence corresponds to a triangulation of a convex $n$-gon by non-intersecting diagonals.

Thus, as defined above, quiddity sequences are the exact same thing as $\eta$-sequences.

\vspace{6 pt}
Given a quiddity ($\eta$-)sequence of positive integers

\vspace{6 pt}
$(a_0,a_1,\dots,a_{n-1}) $

\vspace{6 pt}
\noindent the components of the respective frieze pattern of integers can be easily calculated. Since for the first positive columns $j$, the recursive relation 

\vspace{6 pt}
$ \varphi (k+1,j)=a_k\varphi (k,j)-\varphi (k-1,j) $

\vspace{6 pt}
\noindent is valid (see \cite{Morier-Genoud}, Prop.1.5 or \cite{Leighton}, \S 3) and, as a consequence,

\vspace{6 pt}
$\varphi (k+1,j)=\left( \begin{array}{cccccc} a_j & 1 & 0 & \dots & 0 & 0 \\ 1 & a_{j+1} & 1 & \dots & 0 & 0 \\ 0 & 1 & \ddots & \dots & 0 & 0 \\ \vdots & \vdots & \vdots & \ddots & \vdots  & \vdots \\ 0 & 0 & 0 & \dots & a_{k+j-2} & 1 \\ 0 & 0 & 0 & \dots & 1 & a_{k+j-1} \end{array} \right). $

\vspace{9 pt}
In \cite{Bergeron} (see also \cite{Morier-GenoudOvsienko}) the $SL_2\bf -tilings\rm $ have been defined as double-indexed sequences of integers
$\{\alpha_{i,j}\}_{i\in  \mathbb{Z}, j\in  \mathbb{Z}}$
satisfying the following condition:

\vspace{3pt}
\noindent (FI2) $\alpha_{i,j}\alpha_{i+1,j+1}-\alpha_{i,j+1}\alpha_{i+1,j}=1$  for each $i, j\in  \mathbb{Z}$.

\vspace{3pt}
A tiling is \bf positive \rm if all its components are strictly positive integers

\vspace{3pt}
Notice that a 0 appearing as a component of the tiling must have a $-1$ as neighbor; as a consequence, positive tilings are easy to be described, and their description is based on the following two, easy to prove, Lemmas.

\vspace{6 pt}
\bf Lemma 1\rm . If
\begin{equation*}
\begin{array}{ccc}
\alpha_{i-1,j-1} & \alpha_{i-1,j} & \alpha_{i-1,j+1} \\
\alpha_{i,j-1} & \alpha_{i,j} & \alpha_{i,j+1} \\
\end{array} \left(\rm res.\it 
\begin{array}{cc}
\alpha_{i-1,j-1} & \alpha_{i-1,j} \\
\alpha_{i,j-1} & \alpha_{i,j} \\
\alpha_{i+1,j-1} & \alpha_{i+1,j} \\
\end{array} \right)
\end{equation*}
\noindent are contiguous components of a positive tiling, then there exists a number $k_j\in  \mathbb{N}^*$ (res.  $l_i\in  \mathbb{N}^*$) such that:
\begin{equation*}
\begin{array}{l}
k_j\alpha_{i-1,j}=\alpha_{i-1,j-1}+\alpha_{i-1,j+1} \\
k_j\alpha_{i,j} = \alpha_{i,j-1}+\alpha_{i,j+1} \\
\end{array} \left(\rm res.\it 
\begin{array}{l}
\l_i\alpha_{i,j-1} = \alpha_{i-1,j-1}+ \alpha_{i+1,j-1}\\
\l_i\alpha_{i,j} = \alpha_{i-1,j}+ \alpha_{i+1,j}\\
\end{array} \right). \square
\end{equation*} 
In case $k_j=2$ (res. $l_i=2$), the three consecutive components in row $i$ (res. column $j$), namely $\alpha_{i,j-1},\alpha_{i,j},\alpha_{i,j+1}$ (res. $\alpha_{i-1,j},\alpha_{i,j},\alpha_{i+1,j}$) form an arithmetical progression, and the same is true for the corresponding three components in row $i-1$ (res. column $j-1$). 
In case $k_j\ne 2$ (res. $l_i\ne 2$), the tiling is said to be \bf fractured \rm at column $j$ (res. row $i$). This definition is consistent, due to the following

\vspace{6 pt}
\bf Lemma 2. \rm (See \cite{Morier-GenoudOvsienko}, Proof of Theorem 2.) If
\begin{equation*}
\begin{array}{ccc}
\alpha_{i-1,j-1} & \alpha_{i-1,j} & \alpha_{i-1,j+1} \\
\alpha_{i,j-1} & \alpha_{i,j} & \alpha_{i,j+1} \\
\alpha_{i+1,j-1} & \alpha_{i+1,j} \\
\end{array}
\end{equation*}
\noindent are contiguous components of a positive tiling, and if $k_j, l_i\in  \mathbb{N}^*$  are obtained as in Lemma 1 above, then
\begin{equation*}
\begin{array}{l}
\alpha_{i+1,j+1} = k_j\alpha_{i+1,j}-\alpha_{i+1,j-1} \\
\alpha_{i+1,j+1} = l_i\alpha_{i,j+1}-\alpha_{i-1,j+1}
\end{array}. \ \ \square
\end{equation*}

Hence, given a positive tiling $\{\alpha_{i,j}\}_{i\in  \mathbb{Z}, j\in \mathbb{Z}}$, two infinite vectors $\{ k_j\}_{ j\in \mathbb{Z}}$ and $\{ l_i\}_{ i\in \mathbb{Z}}$ of ``factors" are obtained. If we select only those components that correspond to fractures (at columns, res. at rows), the domain $ \mathbb{Z}\times \mathbb{Z}$ splits into rectangular subdomains (possibly infinite), and in each such subdomain the respective components of the tiling are in arithmetical progressions both on rows and on columns. The tiling is fully determined by these vectors and by an ``initial" $SL_2( \mathbb{Z})$-matrix, for example
\begin{equation*}
\left( \begin{array}{cc}
\alpha_{0,0} & \alpha_{0,1} \\
\alpha_{1,0} & \alpha_{1,1}
\end{array} \right) .
\end{equation*}

Hence, given a positive tiling, there must be at least one fractured column and a fractured row. The following example
\begin{equation*}
\begin{array}{ccccccc}
 & & & \vdots & & & \\
 & 10 & 7 & 4 & 5 & 6 & \\
 & 7 & 5 & 3 & 4 & 5 & \\
 \dots & 4 & 3 & 2 & 3 & 4 & \dots \\
 & 5 & 4 & 3 & 5 & 7 & \\
 & 6 & 5 & 4 & 7 & 10 & \\
 & 5 & 4 & 3 & 5 & 7 & \\
 & & & \vdots & & & \\
\end{array}
\end{equation*}
\noindent shows such a tiling ($\alpha_{i,j}=|i|+|j|+2 \rm \ when\ \it  i\cdot  j<\rm 0,\ and\ \it = |i\cdot  j|+|i|+|j|+\rm 2 \ when\ \it i\cdot j\ge \rm 0$), the fractures being at row=column=0.

\subsection {Frieze patterns of matrices}

\bf Definition 4\rm . By analogy to the definition 2 above, we will call \bf frieze pattern of matrices \rm a double-indexed family
$\{F(i,j)\}_{i\in  \mathbb{N}, j\in  \mathbb{Z}}$
of $2 \times 2$ invertible matrices $F(i,j)$ satisfying two conditions:

\vspace{3pt}
\noindent (FM1) $F(0,j)=F(0,k)$ for all $j,k\in \mathbb{Z}$;

\vspace{3pt} 
\noindent (FM2) (diamond rule) $F(i,j)=F(i-1,j+1)\cdot F(i-2,j+1)^{-1}\cdot F(i-1,j)$

\vspace{3pt}
 for all $i\ge 2, j\in  \mathbb{Z}$.

\vspace{6pt}
Of course, in a frieze pattern of matrices all values in rows 2, 3, ... are well determined provided the 1-st row is known. More precisely, if $M(0,j)=X^{-1}$ and $M(1,j)=A_j$, it is easy to find out that 

\vspace{6pt}
$ M(i,j)=A_{i+j-1}XA_{i+j-2}X\cdot \cdot \cdot XA_j \rm \ \ for \ \it i\ge \rm 2. $

\vspace{6pt}
\noindent And, if all the matrices in rows 0 and 1 are unimodular, then it is clear that all components of such a frieze pattern of matrices belong to $SL_2(\mathbb{Z})$.

\vspace{6pt}
Let us consider the interesting special case $M(0,j)=-S$ for all $j\in \mathbb{Z}$, $M(1,j)=U^{a_j}$ where $a_j\in \mathbb{Z}, j\in \mathbb{Z}$. Then it is clear that

\vspace{6pt}
$ M(i,j)=U^{a_{i+j-1}} S U^{a_{i+j-2}} S \cdot \cdot \cdot  S U^{a_j} \rm \ \ for \  \it i\ge \rm 2,j\in \mathbb{Z} .$

\vspace{6pt}
Denote by $\cal M$  the family of these special freeze patterns of matrices.

\vspace{3pt}
Now, let $\{a_j\}_{j\in  \mathbb{Z}}$ be a sequence of positive integers and consider the frieze pattern (of integers) $\{q(i,j)\}$  generated by putting this sequence as the 2-nd row. By using the diamond rule for $i=2$, we find out that, formally, $q(-1,j)=-1$ for each $j\in \mathbb{Z}$.

Denote by $\{Q(i,j)\}$ the $2\times 2$ integer matrix

\vspace{6pt}
 $\left( \begin{array}{cc} q(i,j) & q(i-1,j+1) \\ q(i+1,j) & q(i,j+1) \end{array} \right) .$

\vspace{6pt}
By the diamond rule {(FI2)}, the matrix $Q(i,j)$ has determinant 1, thus $Q(i,j) \in SL_2(\mathbb{Z})$.

\vspace{6pt}
Formally, $Q(0,j)=-S$ and $Q(1,j)=U^{a_j}$ for all $j$. Applying the diamond rule for matrices, one has $Q(2,j)=U^{a_{j+1}} S U^{a_j}.$ 

Continuing, by repeatedly applying the diamond rule for matrices, one has $Q(i,j)=U^{a_{i+j-1}} S U^{a_{i+j-2}} S \cdot \cdot \cdot  S U^{a_j} \rm \ \ for \  \it i\ge \rm 2 .$

But $q(i,j)$ is recovered as the lower-left component of $Q(i-1,j)$, hence there is a one-to-one correspondence between the family 
$\cal M$ and the sequences $\{a_j\}$ of positive integers.

\vspace{6pt}
Now, if $(a_0,a_1,\dots,a_{n-1})$ is a quiddity sequence, the corresponding special frieze pattern of integers has a row $r>2$  entirely composed of 1's. Formally, the row $r+1$ is composed of 0's, and the next row is composed of $-1$'s. Thus, in the corresponding frieze pattern of matrices one obtains, for the row $r+1$:

\vspace{6pt}
 $F(r+1,j)=\left( \begin{array}{cc} 0 & 1 \\ -1 & 0 \end{array} \right)=S .$

\vspace{6pt}
Conversely, if for a sequence $\{a_j\}$ we construct the corresponding special frieze pattern of matrices and in this freeze pattern the row $r+1$ is composed only of $S$, then from $F(r+1,j)=S$ it follows that

\vspace{6pt}
$U^{a_{r+j}} S U^{a_{r+j-1}} S \cdot \cdot \cdot  S U^{a_{j+1}}SU^{a_j}=S.$

\vspace{6pt}
In particular, $(U^{a_{r}} S U^{a_{r-1}} S \cdot \cdot \cdot  S U^{a_{1}}S)(U^{a_0}S)=(U^{a_{r+1}} S )(U^{a_{r}} S \cdot \cdot \cdot  S U^{a_{1}})$; hence, from Lemma 1 we obtain $a_{r+1}=a_0$. In consequence, $(a_0,a_1,\dots,a_{r})$ is a quiddity sequence. This proves the following

\vspace{3pt}
\bf Theorem 1\rm . Under the notations above, an expression $U^{a_{n-1}} S U^{a_{n-2}} S \cdot\cdot\cdot  S U^{a_0}$ reduces to $-I$ if and only if  $(a_0,a_1,\dots,a_{n-1})$ is an $\eta$-sequence.  $\square$

  \section{Embedding in quiddity sequences}


In what follows we study the possibility of embedding a finite sequence in a quiddity sequence. 
First, any integer $A$ can be embedded in an $\eta$-sequence; this is obvious, because it is enough to apply repeatedly Rule 3 to the initial $\eta$-sequence $(1,1,1)$ resulting in the (left-sided) \bf quiddity fan sequence \rm $(A,1,\underbrace{2,2,\dots ,2}_{\it A-\rm 1\ \ copies},1)$. For completeness, a (right-sided)  \bf quiddity fan sequence \rm is of the form $(1,\underbrace{2,2,\dots ,2}_{\it A-\rm 1\ \ copies},1,A)$.

Alternatively, a sequence $1,A$ can be supplemented by the sequence $1,2,\dots ,2$ ($A-1$ copies of 2) to form a quiddity sequence.

\vspace{6pt}
Let us consider now sequences $1,A_1,\dots ,A_n$ with $n\ge 1$ and all the components $A_1,\dots,A_n\ge 2$, which will be referred to as \bf basic\rm.

\vspace{3pt}
\bf Proposition 2\rm . For each basic sequence $\alpha$ there exists a ``supplementary" basic sequence $\overline\alpha$ such that the concatenation $\alpha \overline\alpha$ is a quiddity sequence.

\vspace{3pt}
\it Proof\rm . Decompose $\alpha$ into consecutive subsequences, by separating the components 2 from the others, as follows:
$$\alpha\ =\ (1,\underbrace{2,\dots ,2}_{x_0\rm \ copies},A_1,\underbrace{2,\dots ,2}_{x_1\rm \ copies},A_2,\dots ,\underbrace{2,\dots ,2}_{x_{k-1}\rm \ copies},A_k,\underbrace{2,\dots ,2}_{x_k\rm \ copies})$$
\noindent where $A_1,A_2,\dots ,A_k\ge 3$ and $x_0,x_1,\dots,x_k\ge 0$.

\vspace{6pt}
Let us distinguish several cases.

The case $k=1,x_0=x_1=0$ has been treated above. The supplement $\overline\alpha$ is $(1,2,\dots ,2)$ with $A-1$ copies of 2.

Dually, the case $k=0,x_0\ge 1$, the supplement  $\overline\alpha$ is $(1,x_0+1)$.

For the other cases, the supplement will be constructed step by step, taking into account all the subsequences of $\alpha$, but in reverse order. The subsequences $\underbrace{2,\dots ,2}_{x_l\rm \ copies}$ will produce a singleton subsequence $X_l=x_l+3$, and the subsequences $A_l$ will produce subsequences (possible void) $\underbrace{2,\dots ,2}_{A_l-3\rm \ copies}$. There are four exceptions:

1) in case $x_k\ne 0$, the second component in $\overline\alpha$ (after the initial 1) will be $X_k=x_k+2$;

2) in case $x_k=0$, the $A_k$ will produce the subsequence $\underbrace{2,\dots ,2}_{A_k-2\rm \ copies}$;

3) in case $x_0\ne 0$, the last term in $\overline\alpha$ will be $X_0=x_0+3$, and

4) in case $x_0=0$, $A_1$ will produce  the subsequence $\underbrace{2,\dots ,2}_{A_1-2\rm \ copies}$.

Thus, the supplement of $\alpha$ is  

\noindent$\overline\alpha \ =\ (1,X_k,2,... ,2,X_{k-1},...,\underbrace{2,...,2}_{A_l-3\rm \ copies},X_{l-1},...,2,... ,2,X_0,2,...,2).\square$


Obviously, the supplements are involutive: $\overline{\overline\alpha }=\alpha$.

For example, the basic sequence $\overline\alpha$ (1,6,3,2,4,2,2,2,4) is the supplement of the basic sequence $\alpha$ (1,2,2,6,2,4,3,2,2,2,2). In the dual representation, the elements of a basic sequence $\alpha$ are given by the counts of consecutive internal nodes in the pre-order step (left side of tree) of the binary tree traversal. The supplementary sequence $\overline\alpha$ is similar except that its elements are the counts of internal nodes in the post-order step (right side of tree) of the traversal.

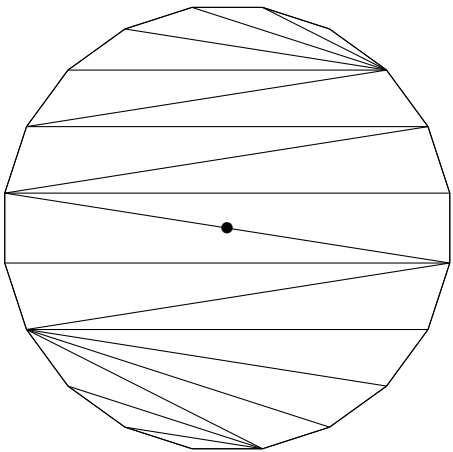
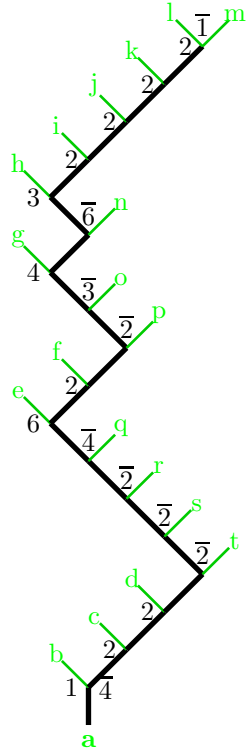
\begin{tikzpicture}

\polygon[scale=3,xshift=2cm,yshift=4cm]{20}{{11/13,10/13,9/13,8/13,7/13,7/14,6/14,6/15,6/16,5/16,4/16,4/17,4/18,4/19,4/20,3/20,2/20}}

\begin{scope}[scale=.5,xshift=14cm,yshift=-3cm]
\draw (10,20)		node[anchor=north] {\small\color{green}\bf a};
\draw [line width=2pt] (10,20)-- (10,21);
\draw [color=zzttgg, line width=1pt] (10,21)-- (9.3,21.7);
\draw (9.15,21.4)		node[anchor=south] {\small\color{green}b};
\draw (10,21)		node[anchor=east] {\small1};
\draw [line width=2pt] (10,21)-- (11,22);
\draw [color=zzttgg, line width=1pt] (11,22)-- (10.3,22.7);
\draw (10.15,22.4)		node[anchor=south] {\small\color{green}c};
\draw (11,22)		node[anchor=east] {\small2};
\draw [line width=2pt] (11,22)-- (12,23);
\draw [color=zzttgg, line width=1pt] (12,23)-- (11.3,23.7);
\draw (11.15,23.4)		node[anchor=south] {\small\color{green}d};
\draw (12,23)		node[anchor=east] {\small2};
\draw [line width=2pt] (12,23)-- (13,24);
\draw [line width=2pt] (13,24)-- (12,25);
\draw [line width=2pt] (12,25)-- (11,26);
\draw [line width=2pt] (11,26)-- (10,27);
\draw [line width=2pt] (10,27)-- (9,28);
\draw [color=zzttgg, line width=1pt] (9,28)-- (8.3,28.7);
\draw (8.15,28.4)		node[anchor=south] {\small\color{green}e};
\draw (9,28)		node[anchor=east] {\small6};
\draw [line width=2pt] (9,28)-- (10,29);
\draw [color=zzttgg, line width=1pt] (10,29)-- (9.3,29.7);
\draw (9.15,29.4)		node[anchor=south] {\small\color{green}f};
\draw (10,29)		node[anchor=east] {\small2};
\draw [line width=2pt] (10,29)-- (11,30);
\draw [line width=2pt] (11,30)-- (10,31);
\draw [line width=2pt] (10,31)-- (9,32);
\draw [color=zzttgg, line width=1pt] (9,32)-- (8.3,32.7);
\draw (8.15,32.4)		node[anchor=south] {\small\color{green}g};
\draw (9,32)		node[anchor=east] {\small4};
\draw [line width=2pt] (9,32)-- (10,33);
\draw [line width=2pt] (10,33)-- (9,34);
\draw [color=zzttgg, line width=1pt] (9,34)-- (8.3,34.7);
\draw (8.15,34.4)		node[anchor=south] {\small\color{green}h};
\draw (9,34)		node[anchor=east] {\small3};
\draw [line width=2pt] (9,34)-- (10,35);
\draw [color=zzttgg, line width=1pt] (10,35)-- (9.3,35.7);
\draw (9.15,35.4)		node[anchor=south] {\small\color{green}i};
\draw (10,35)		node[anchor=east] {\small2};
\draw [line width=2pt] (10,35)-- (11,36);
\draw [color=zzttgg, line width=1pt] (11,36)-- (10.3,36.7);
\draw (10.15,36.4)		node[anchor=south] {\small\color{green}j};
\draw (11,36)		node[anchor=east] {\small2};
\draw [line width=2pt] (11,36)-- (12,37);
\draw [color=zzttgg, line width=1pt] (12,37)-- (11.3,37.7);
\draw (11.15,37.4)		node[anchor=south] {\small\color{green}k};
\draw (12,37)		node[anchor=east] {\small2};
\draw [line width=2pt] (12,37)-- (13,38);
\draw [color=zzttgg, line width=1pt] (13,38)-- (12.3,38.7);
\draw (12.15,38.4)		node[anchor=south] {\small\color{green}l};
\draw (13,38)		node[anchor=east] {\small2};
\draw [color=zzttgg, line width=1pt] (13,38)-- (13.7,38.7);
\draw (13.85,38.4)		node[anchor=south] {\small\color{green}m};
\draw (13,38)		node[anchor=south] {\small$\overline{1}$};
\draw [color=zzttgg, line width=1pt] (10,33)-- (10.7,33.7);
\draw (10.85,33.4)		node[anchor=south] {\small\color{green}n};
\draw (10,33)		node[anchor=south] {\small$\overline{6}$};
\draw [color=zzttgg, line width=1pt] (10,31)-- (10.7,31.7);
\draw (10.85,31.4)		node[anchor=south] {\small\color{green}o};
\draw (10,31)		node[anchor=south] {\small$\overline{3}$};
\draw [color=zzttgg, line width=1pt] (11,30)-- (11.7,30.7);
\draw (11.85,30.4)		node[anchor=south] {\small\color{green}p};
\draw (11,30)		node[anchor=south] {\small$\overline{2}$};
\draw [color=zzttgg, line width=1pt] (10,27)-- (10.7,27.7);
\draw (10.85,27.4)		node[anchor=south] {\small\color{green}q};
\draw (10,27)		node[anchor=south] {\small$\overline{4}$};
\draw [color=zzttgg, line width=1pt] (11,26)-- (11.7,26.7);
\draw (11.85,26.4)		node[anchor=south] {\small\color{green}r};
\draw (11,26)		node[anchor=south] {\small$\overline{2}$};
\draw [color=zzttgg, line width=1pt] (12,25)-- (12.7,25.7);
\draw (12.85,25.4)		node[anchor=south] {\small\color{green}s};
\draw (12,25)		node[anchor=south] {\small$\overline{2}$};
\draw [color=zzttgg, line width=1pt] (13,24)-- (13.7,24.7);
\draw (13.85,24.4)		node[anchor=south] {\small\color{green}t};
\draw (13,24)		node[anchor=south] {\small$\overline{2}$};
\draw (10,21)		node[anchor=west] {\small$\overline{4}$};

\end{scope}

\end{tikzpicture}\\

\vspace{6pt}
If $\alpha$ is not basic, the result above may not be valid. Particularly when $(2,1,2)$ is contained as a subsequence in $\alpha$ (which is different from the repeated $1,2$), it is clear that $\alpha$ cannot be inserted in a quiddity sequence. However, we can extend the Proposition 2 above to sequences containing more than a $1$ inside. Namely to ``super-basic" sequences.

\vspace{6pt}
\bf Definition 5\rm . A sequence $\alpha\ =\ (1,A_1,\dots,A_n)$ is considered \bf super-basic \rm if it is basic, $n>1$ and $A_1,A_n>2$.

\vspace{6pt}
Let $\alpha =(1,A_1,\dots,A_n)$ and $\beta =(1,B_1,\dots,B_m)$ two super-basic sequences. Then $\hat\alpha =(1,A_1,\dots,A_{n-1},A_n-1)$ and $\hat\beta =(1,B_1-1,B_2,\dots,B_m)$ are basic. Consider their respective supplements $\xi\ =\ (1,X_1,\dots,X_p)$ for $\hat\alpha$, $\eta\ =\ (1,Y_1,\dots,Y_q)$ for $\hat\beta$. 

Since $A_1>2,A_n-1\ge 2$ and $B_1-1\ge 2, B_m>2$, the sequence

\vspace{6 pt}
\hspace{40 pt}$\gamma\ = \ (1,A_1,\dots,A_n-1,B_1-1,\dots,B_m)$

\vspace{6 pt}
\noindent is (super)basic and its supplement is

\vspace{6 pt}
\hspace{50 pt}$\zeta\ =\ (1,Y_1,\dots,Y_{q-1},Y_q+X_1-1,X_2,\dots,X_p).$

\vspace{6 pt}
Now, applying the (expansion) Rule 3 to the concatenated sequence $\gamma\zeta$, at the pair $(A_n-1,B_1-1)$, it follows that the concatenation $\alpha\beta$ extends to the quiddity sequence $\alpha\beta\zeta$.

\vspace{6pt}
It is easy to extend the reasoning above to prove the following:

\bf Proposition 3\rm . If $\alpha_1,\alpha_2,...,\alpha_s$ are super-basic sequences, then their concatenation $\alpha_1\alpha_2\dots\alpha_s$ can be extended to a quiddity sequence. $\square$

  \section{Similarity type of frieze patterns}

Let us look at the following three examples of frieze patterns (of numbers):

\noindent$...\ 1\ \ 1\ \ 1\ \ 1\ \ 1\ \ 1\ \ 1\ ...\ \ \ ...\ 1\ \ 1\ \ 1\ \ 1\ \ 1\ \ 1\ \ 1\ ...\ \ \ \ ...\ 1\ \ 1\ \ 1\ \ 1\ \ 1\ \ 1\ \ 1\ ...$

\noindent$\ \ ...\ 4\ \ 2\ \ 1\ \ 3\ \ 2\ \ 2\ \ 1\ ...\ \ \ ...\ 1\ \ 2\ \ 2\ \ 3\ \ 1\ \ 2\ \ 4\ ...\ \ \ \ ...\ 4\ \ 1\ \ 3\ \ 1\ \ 3\ \ 2\ \ 1\ ...$

\noindent$...\ 3\ \ 7\ \ 1\ \ 2\ \ 5\ \ 3\ \ 1\ ...\ \ \ ...\ 3\ \ 1\ \ 3\ \ 5\ \ 2\ \ 1\ \ 7\ ...\ \ \ \ ...\ 3\ \ 3\ \ 2\ \ 2\ \ 2\ \ 5\ \ 1\ ...$

\noindent$\ \ ...\ 5\ \ 3\ \ 1\ \ 3\ \ 7\ \ 1\ \ 2\ ...\ \ \ ...\ 2\ \ 1\ \ 7\ \ 3\ \ 1\ \ 3\ \ 5\ ...\ \ \ \ ...\ 2\ \ 5\ \ 1\ \ 3\ \ 3\ \ 2\ \ 2\ ...$

\noindent$...\ 3\ \ 2\ \ 2\ \ 1\ \ 4\ \ 2\ \ 1\ ...\ \ \ ...\ 3\ \ 1\ \ 2\ \ 4\ \ 1\ \ 2\ \ 2\ ...\ \ \ \ ...\ 1\ \ 3\ \ 2\ \ 1\ \ 4\ \ 1\ \ 3\ ...$

\noindent$\ \ ...\ 1\ \ 1\ \ 1\ \ 1\ \ 1\ \ 1\ \ 1\ ...\ \ \ ...\ 1\ \ 1\ \ 1\ \ 1\ \ 1\ \ 1\ \ 1\ ...\ \ \ \ ...\ 1\ \ 1\ \ 1\ \ 1\ \ 1\ \ 1\ \ 1\ ...$

\noindent$\hspace{50 pt} \it (I)\hspace{112 pt} (I') \hspace {114 pt} (II) \rm $

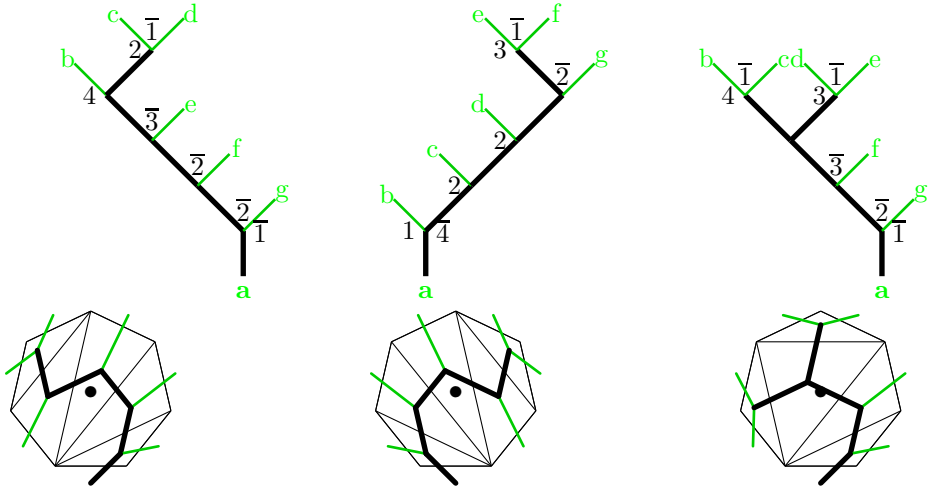
\vspace{6 pt}

\begin{tikzpicture}[scale=.6]
\begin{scope}
\draw (3.5,7)		node[anchor=north] {\small\color{green}\bf a};
\draw [line width=2pt] (3.5,7)-- (3.5,8);
\draw [line width=2pt] (3.5,8)-- (2.5,9);
\draw [line width=2pt] (2.5,9)-- (1.5,10);
\draw [line width=2pt] (1.5,10)-- (0.5,11);
\draw [color=zzttgg, line width=1pt] (0.5,11)-- (-0.2,11.7);
\draw (-0.35,11.4)		node[anchor=south] {\small\color{green}b};
\draw (0.5,11)		node[anchor=east] {\small4};
\draw [line width=2pt] (0.5,11)-- (1.5,12);
\draw [color=zzttgg, line width=1pt] (1.5,12)-- (0.8,12.7);
\draw (0.65,12.4)		node[anchor=south] {\small\color{green}c};
\draw (1.5,12)		node[anchor=east] {\small2};
\draw [color=zzttgg, line width=1pt] (1.5,12)-- (2.2,12.7);
\draw (2.35,12.4)		node[anchor=south] {\small\color{green}d};
\draw (1.5,12)		node[anchor=south] {\small$\overline{1}$};
\draw [color=zzttgg, line width=1pt] (1.5,10)-- (2.2,10.7);
\draw (2.35,10.4)		node[anchor=south] {\small\color{green}e};
\draw (1.5,10)		node[anchor=south] {\small$\overline{3}$};
\draw [color=zzttgg, line width=1pt] (2.5,9)-- (3.2,9.7);
\draw (3.35,9.4)		node[anchor=south] {\small\color{green}f};
\draw (2.5,9)		node[anchor=south] {\small$\overline{2}$};
\draw [color=zzttgg, line width=1pt] (3.5,8)-- (4.2,8.7);
\draw (4.35,8.4)		node[anchor=south] {\small\color{green}g};
\draw (3.5,8)		node[anchor=south] {\small$\overline{2}$};
\draw (3.5,8)		node[anchor=west] {\small$\overline{1}$};
\end{scope}

\begin{scope}[xshift=4cm]

\draw (3.5,7)		node[anchor=north] {\small\color{green}\bf a};
\draw [line width=2pt] (3.5,7)-- (3.5,8);
\draw [color=zzttgg, line width=1pt] (3.5,8)-- (2.8,8.7);
\draw (2.65,8.4)		node[anchor=south] {\small\color{green}b};
\draw (3.5,8)		node[anchor=east] {\small1};
\draw [line width=2pt] (3.5,8)-- (4.5,9);
\draw [color=zzttgg, line width=1pt] (4.5,9)-- (3.8,9.7);
\draw (3.65,9.4)		node[anchor=south] {\small\color{green}c};
\draw (4.5,9)		node[anchor=east] {\small2};
\draw [line width=2pt] (4.5,9)-- (5.5,10);
\draw [color=zzttgg, line width=1pt] (5.5,10)-- (4.8,10.7);
\draw (4.65,10.4)		node[anchor=south] {\small\color{green}d};
\draw (5.5,10)		node[anchor=east] {\small2};
\draw [line width=2pt] (5.5,10)-- (6.5,11);
\draw [line width=2pt] (6.5,11)-- (5.5,12);
\draw [color=zzttgg, line width=1pt] (5.5,12)-- (4.8,12.7);
\draw (4.65,12.4)		node[anchor=south] {\small\color{green}e};
\draw (5.5,12)		node[anchor=east] {\small3};
\draw [color=zzttgg, line width=1pt] (5.5,12)-- (6.2,12.7);
\draw (6.35,12.4)		node[anchor=south] {\small\color{green}f};
\draw (5.5,12)		node[anchor=south] {\small$\overline{1}$};
\draw [color=zzttgg, line width=1pt] (6.5,11)-- (7.2,11.7);
\draw (7.35,11.4)		node[anchor=south] {\small\color{green}g};
\draw (6.5,11)		node[anchor=south] {\small$\overline{2}$};
\draw (3.5,8)		node[anchor=west] {\small$\overline{4}$};

\end{scope}

\begin{scope}[xshift=14cm]
\draw (3.5,7)		node[anchor=north] {\small\color{green}\bf a};
\draw [line width=2pt] (3.5,7)-- (3.5,8);
\draw [line width=2pt] (3.5,8)-- (2.5,9);
\draw [line width=2pt] (2.5,9)-- (1.5,10);
\draw [line width=2pt] (1.5,10)-- (0.5,11);
\draw [color=zzttgg, line width=1pt] (0.5,11)-- (-0.2,11.7);
\draw (-0.35,11.4)		node[anchor=south] {\small\color{green}b};
\draw (0.5,11)		node[anchor=east] {\small4};
\draw [color=zzttgg, line width=1pt] (0.5,11)-- (1.2,11.7);
\draw (1.35,11.4)		node[anchor=south] {\small\color{green}c};
\draw (0.5,11)		node[anchor=south] {\small$\overline{1}$};
\draw [line width=2pt] (1.5,10)-- (2.5,11);
\draw [color=zzttgg, line width=1pt] (2.5,11)-- (1.8,11.7);
\draw (1.65,11.4)		node[anchor=south] {\small\color{green}d};
\draw (2.5,11)		node[anchor=east] {\small3};
\draw [color=zzttgg, line width=1pt] (2.5,11)-- (3.2,11.7);
\draw (3.35,11.4)		node[anchor=south] {\small\color{green}e};
\draw (2.5,11)		node[anchor=south] {\small$\overline{1}$};
\draw [color=zzttgg, line width=1pt] (2.5,9)-- (3.2,9.7);
\draw (3.35,9.4)		node[anchor=south] {\small\color{green}f};
\draw (2.5,9)		node[anchor=south] {\small$\overline{3}$};
\draw [color=zzttgg, line width=1pt] (3.5,8)-- (4.2,8.7);
\draw (4.35,8.4)		node[anchor=south] {\small\color{green}g};
\draw (3.5,8)		node[anchor=south] {\small$\overline{2}$};
\draw (3.5,8)		node[anchor=west] {\small$\overline{1}$};
\end{scope}

\end{tikzpicture}\\
\begin{tikzpicture}[scale=.6]
\begin{scope}[xshift=0cm]
\polygon[scale=2,xshift=5cm,yshift=5cm]{7}{{2/4,1/4,1/5,1/6}}
\draw [color=zzttgg, line cap=round, line join=round, line width=1pt] (8.829,10.934)-- (8.148,10.423);
\draw [color=zzttgg, line cap=round, line join=round, line width=1pt] (8.829,10.934)-- (9.176,11.712);
\draw [color=zzttgg, line cap=round, line join=round, line width=1pt] (9.061,9.918)-- (8.515,8.815);
\draw [line cap=round, line join=round, line width=2pt] (9.061,9.918)-- (8.829,10.934);
\draw [line cap=round, line join=round, line width=2pt] (10.232,10.482)-- (9.061,9.918);
\draw [color=zzttgg, line cap=round, line join=round, line width=1pt] (10.232,10.482)-- (10.824,11.712);
\draw [line cap=round, line join=round, line width=2pt] (10.882,9.667)-- (10.232,10.482);
\draw [color=zzttgg, line cap=round, line join=round, line width=1pt] (10.882,9.667)-- (11.852,10.423);
\draw [line cap=round, line join=round, line width=2pt] (10.650,8.650)-- (10.882,9.667);
\draw [color=zzttgg, line cap=round, line join=round, line width=1pt] (10.650,8.650)-- (11.485,8.815);
\draw [line cap=round, line join=round, line width=2pt] (10.650,8.650)-- (10.000,8.000);
\end{scope}

\begin{scope}[xshift=8cm]

\polygon[scale=2,xshift=5cm,yshift=5cm]{7}{{4/6,4/7,3/7,2/7}}
\draw [color=zzttgg, line cap=round, line join=round, line width=1pt] (11.171,10.934)-- (10.824,11.712);
\draw [color=zzttgg, line cap=round, line join=round, line width=1pt] (11.171,10.934)-- (11.852,10.423);
\draw [line cap=round, line join=round, line width=2pt] (10.939,9.918)-- (11.171,10.934);
\draw [color=zzttgg, line cap=round, line join=round, line width=1pt] (10.939,9.918)-- (11.485,8.815);
\draw [color=zzttgg, line cap=round, line join=round, line width=1pt] (9.768,10.482)-- (9.176,11.712);
\draw [line cap=round, line join=round, line width=2pt] (9.768,10.482)-- (10.939,9.918);
\draw [color=zzttgg, line cap=round, line join=round, line width=1pt] (9.118,9.667)-- (8.148,10.423);
\draw [line cap=round, line join=round, line width=2pt] (9.118,9.667)-- (9.768,10.482);
\draw [color=zzttgg, line cap=round, line join=round, line width=1pt] (9.350,8.650)-- (8.515,8.815);
\draw [line cap=round, line join=round, line width=2pt] (9.350,8.650)-- (9.118,9.667);
\draw [line cap=round, line join=round, line width=2pt] (9.350,8.650)-- (10.000,8.000);
\end{scope}

\begin{scope}[xshift=16cm]
\polygon[scale=2,xshift=5cm,yshift=5cm]{7}{{1/3,3/5,1/5,1/6}}
\draw [color=zzttgg, line cap=round, line join=round, line width=1pt] (8.540,9.667)-- (8.515,8.815);
\draw [color=zzttgg, line cap=round, line join=round, line width=1pt] (8.540,9.667)-- (8.148,10.423);
\draw [color=zzttgg, line cap=round, line join=round, line width=1pt] (10.000,11.498)-- (9.176,11.712);
\draw [color=zzttgg, line cap=round, line join=round, line width=1pt] (10.000,11.498)-- (10.824,11.712);
\draw [line cap=round, line join=round, line width=2pt] (9.711,10.231)-- (8.540,9.667);
\draw [line cap=round, line join=round, line width=2pt] (9.711,10.231)-- (10.000,11.498);
\draw [line cap=round, line join=round, line width=2pt] (10.882,9.667)-- (9.711,10.231);
\draw [color=zzttgg, line cap=round, line join=round, line width=1pt] (10.882,9.667)-- (11.852,10.423);
\draw [line cap=round, line join=round, line width=2pt] (10.650,8.650)-- (10.882,9.667);
\draw [color=zzttgg, line cap=round, line join=round, line width=1pt] (10.650,8.650)-- (11.485,8.815);
\draw [line cap=round, line join=round, line width=2pt] (10.650,8.650)-- (10.000,8.000);
\end{scope}

\end{tikzpicture}\\

The similarity between $\it (I)\rm$  and $\it (I')\rm$  is obvious but $\it (II)\rm$ is clearly a different type. In what follows we study this type of similarity, being well aware that frieze patterns are determined by quiddity sequences.

\vspace{6pt}
The cyclic permutation $\gamma (i)=i+1 (\rm modulo\ \it n)$ and the symmetry $\sigma (i)=n-i$ generate a dihedral group $\mathbb{D}_n$, which acts obviously over the set of quiddity sequences of length $n$. The orbits of this action (see \cite{Fraleigh}, pag. 158) will be called \bf similarity types \rm (of quiddity sequences, thus of frieze patterns) and we will denote by $K_n$ the number of these orbits.

We already know that quiddity sequences $(a_i,a_{i+1},...,a_{i-1})$ obtained from $(a_0,a_1,...,a_{n-1})$ by cyclic permutation of indices determine the same frieze pattern, thus they can be considered as being similar. 
In addition, two quiddity sequences $\bf a\it =(a_{\rm 0\it},a_{\rm 1\it},...,a_{n-\rm 1})$ and $\bf b\it =(b_{\rm 0\it},b_{\rm 1\it},...,b_{n-\rm 1})$ (of the same length) are called \bf similar \rm if $\bf b\rm $ is obtained from $\bf a\rm $ as a result of the action of $\sigma$. Notation: $\bf b\it =\overline{\bf a\rm }$.

\vspace{6pt}
A quiddity sequence is called \bf symmetric \rm in case $\overline{\bf a\it }=\bf a\rm$.

Since $\overline{\overline{\bf a\it }}=\bf a\rm $ , the non-symmetric quiddity sequences are grouped in pairs $\{ \bf a\it , \overline{\bf a\rm }\} $.

\vspace{6pt}
Denote by $\mathcal{T}_n$ the set of quiddity sequences of length $n$, and by $\mathcal{S}_n$ its subset composed of symmetric quiddity sequences. Point (i) in the following proposition is obvious.

\vspace{6pt}
\bf Proposition 4\rm . (i) $T_n=2A_n+S_n$, where $T_n=|\mathcal{T}_n|$ and $S_n=|\mathcal{S}_n|$;

(ii) $T_n=\rm C_{\it n -\rm  2}$, where $\rm C_{\it k}$ are the Catalan numbers;

(iii) $S_n=0$ if $n=2m$, and $=\rm C_{\it m \rm -1}$ if $n=2m+1$.

\vspace{3pt}
\it Proof\rm . (ii) is well-known (see [4]).

(iii) For $n=2m+1$, it is obvious that any symmetric quiddity sequence $(a_0,a_1,...,a_{2m})$ is perfectly determined by the ``first half" $(a_0,a_1,...,a_m)$. For $n=2m$, the proof is easy if we adopt the language of triangulations, considering that $a_i$ is the number of triangles supported by the vertex $i$ of a convex $n$-gon. Namely, the segment $\{m-1,m\}$ belongs to a certain triangle $t=\{m-1,m,k\}$ where $k<m-1$ or $k>m$. Its symmetrical is $\overline{t}=\{m-1,m,m-k\}$, thus $\overline{\tau}=\tau$ is impossible. $\square$

\vspace{3pt}
In what follows, a \bf perfect tri-partition \rm  of $n$ is a triple $\{ i,j,k\} $ such that $i\ge j\ge k$ and $n=i+j+k$, and 

(a) if $n=2m$, then either $i=j=m, k=0$, or $i,j<m, k\ge 2$;

(b) if $n=2m+1$, then either $i=j=m,k=1$, or $i,j,k<m$.

\vspace{3pt}
Denote by $\mathcal{P}_n$ the set of perfect tri-partitions of $n$. For an element $\{ i,j,k\}\in\mathcal{P}_n$, six cases are possible:

(A) $n=2m, i=j=m, k=0$.

(B) $n=2m+1, i=j=m, k=1$.

(C) $i>j>k$.

(D) $i=j>k$.

(E) $i>j=k$.

(F) $i=j=k$.

\vspace{3 pt}
Now return to the triangulations of a convex regular $n$-gon. 

The similarity types (of triangulations, thus of frieze patterns) attached to $n$ can be obtained for each of these tri-partitions, according to the asociated case.

\vspace{3pt}
Let us select and fix a vertex. 

In case (A), the $n$-gon is split by its center and this fixed vertex into two (symmetric) $(m+1)$-gons. Thus, any triangulation $\tau\in\mathcal{T}_n$ of the initial $n$-gon  is determined by an ordered pair $(\tau_1,\tau_2)$ of triangulations of these two $(m+1)$-gons. 
To obtain the number of similarity types, $N(m+1,m+1,0)$, in this case, notice that the group $\mathbb{Z}_2\times \mathbb{Z}_2$ acts by symmetries over the pairs $(\tau_1,\tau_2)$. The possible orbits may contain 1, 2 or 4 elements. Taking into account the relation $T_{m+1}=2A_{m+1}+S_{m+1}$, the following formula is valid:

\vspace{3 pt}
$N(m,m,0)= A_{m+1}( A_{m+1}+1)+ A_{m+1}S_{m+1}+\frac{S_{m+1}(S_{m+1}+1)}{2}$.

\vspace{3 pt}
In the other five cases (B-F), the splitting involves three polygons circonding the central triangle, namely a $(i+1)$-gon, a $(j+1)$-gon and a $(k+1)$-gon (except in case (B), when the $(k+1)$-gon is degenerate). 

\vspace{3pt}
To obtain the number of similarity types in case (B), notice that  group $\mathbb{Z}_2$ is acting over the pairs $(\tau_1,\tau_2)$, determining orbits of 1 or 2 elements. Overall, the number of orbits is:

\vspace{3 pt}
$N(m,m,1)=\frac{T_{m+1}(T_{m+1}+1)}{2}$.

\vspace{3 pt}
In case (C), the initial $n$-gon is split into three different-sized polygons (plus the central triangle). Each possible triple $(\tau_i,\tau_j, \tau_k)$, where $\tau_i$ is a triangulation of the resulting $(i+1)$-gon, etc., determines an overall triangularization, of a specific type. It is rather obvious that:

\vspace{3 pt}
$N(i,j,k)= T_{i+1} T_{j+1}T_{k+1}$.

\vspace{3 pt}
In cases (D) and (E) the splitting involves two equal-sized polygons, and a third of different size. To treat both as a single case, suppose an $x$-gon and two $y$-gons are involved, and take into account the relations  $T_x=2A_x+S_x$ and $T_y=2A_y+S_y$. Again,  group $\mathbb{Z}_2$ is acting over the triples $(\tau_x,\tau_{y1},\tau_{y2})$, determining orbits of 1 or 2 elements. It is easy to count the 1-element orbits, they correspond to triples in which $\tau_x$ is symmetrical, $\tau_{y1}$ is arbitrary and $\tau_{y2}=\overline{\tau_{y1}}$. Thus

\vspace{3 pt}
$ N(i,i,k)=\frac{T_{i+1}^2T_{k+1}+T_{i+1}S_{k+1}}{2},\ \  N(i,k,k)=\frac{T_{i+1}T_{k+1}^2+S_{i+1}T_{k+1}}{2}. $

\vspace{3 pt}
Finally, in case (F) the splitting involves three equal-sized polygons (plus the central triangle). The dihedral group $\mathbb{D}_3$ acts over the set of triples  $(\tau',\tau'',\tau''')$, where $\tau',\tau'',\tau''' $ are triagularizations of respective $(i+1)$-gons. The possible orbits may have 1, 2, 3 or 6 elements. Taking into account the decomposition $T_{i+1}=2A_{i+1}+S_{i+1}$, a thorough counting of orbits gives:

\vspace{3 pt}
$N(i,i,i)=\frac{T_{i+1}(T_{i+1}+1)(T_{i+1}+2)}{6}-T_{i+1}A_{i+1}$.

\vspace{3 pt}
We have thus obtained the following.

\bf Proposition 5\rm . The following formula can be used to compute the number of different similarity types of frieze patterns of numbers 

\vspace{3pt}
$ K_n=\sum_{(i,j,k)\in\mathcal{P}_n}N(i,j,k).\square $

\vspace{6pt}
For example, there are 5 perfect tri-partitions of 13, namely 6+6+1, 6+5+2, 6+4+3, 5+5+3, 5+4+4. The corresponding numbers (obtained using the formulas for cases B, C, C, D res. E) are as follows: 903, 588, 420, 196 res. 175. Summing up these terms, one obtains $K_{13}=2282$.

\vspace{6pt}
The following list contains a few values of these numbers:

\noindent \begin{tabular}{|c|rrrrrrrrrrrc|}
\hline
$n$&3&4&5&6&7&8&9&10&11&12&13&...\\
\hline
$T_n$&1&2&5&14&42&132&429&1430&4862&16796&58786&...\\
$S_n$&1&0&1&0&2&0&5&0&14&0&42&...\\
$A_n$&0&1&2&7&20&66&221&715&2424&8398&29372&...\\
\hline
$K_n$&1&1&1&3&4&12&27&82&228&733&2282&...\\
\hline
\end{tabular}

\vspace{6pt}
To give a last example, let us specify that $K_7=4$. Two of the four types of non-similar frieze patterns based on quiddity sequences of length 7 are presented in Figures $(I)$ and $(II)$ above. In fact, these two frieze patterns represent the non-symmetrical ones. The other two (which are symmetrical), are the following:\\

$...\ 1\ \ 1\ \ 1\ \ 1\ \ 1\ \ 1\ \ 1\ ...\ \ \ \ ...\ 1\ \ 1\ \ 1\ \ 1\ \ 1\ \ 1\ \ 1\ ...$

$\ \ ...\ 5\ \ 1\ \ 2\ \ 2\ \ 2\ \ 2\ \ 1\ ...\ \ \ \ ...\ 3\ \ 2\ \ 1\ \ 3\ \ 3\ \ 1\ \ 2\ ...$

$...\ 4\ \ 4\ \ 1\ \ 3\ \ 3\ \ 3\ \ 1\ ...\ \ \ \ ...\ 5\ \ 5\ \ 1\ \ 2\ \ 8\ \ 2\ \ 1\ ...$

$\ \ ...\ 3\ \ 3\ \ 1\ \ 4\ \ 4\ \ 1\ \ 3\ ...\ \ \ \ ...\ 8\ \ 2\ \ 1\ \ 5\ \ 5\ \ 1\ \ 2\ ...$

$...\ 2\ \ 2\ \ 2\ \ 1\ \ 5\ \ 1\ \ 2\ ...\ \ \ \ ...\ 3\ \ 3\ \ 1\ \ 2\ \ 3\ \ 2\ \ 1\ ...$

$\ \ ...\ 1\ \ 1\ \ 1\ \ 1\ \ 1\ \ 1\ \ 1\ ...\ \ \ \ ...\ 1\ \ 1\ \ 1\ \ 1\ \ 1\ \ 1\ \ 1\ ... $ \\

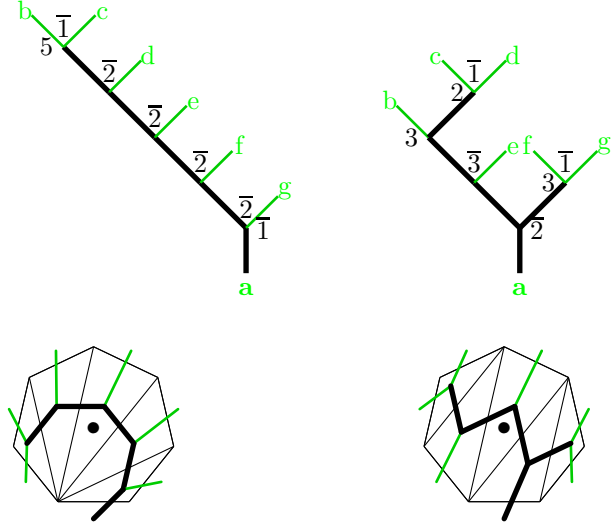
\begin{tikzpicture}[scale=.6]

\begin{scope}
\draw (3.5,7)		node[anchor=north] {\small\color{green}\bf a};
\draw [line width=2pt] (3.5,7)-- (3.5,8);
\draw [line width=2pt] (3.5,8)-- (2.5,9);
\draw [line width=2pt] (2.5,9)-- (1.5,10);
\draw [line width=2pt] (1.5,10)-- (0.5,11);
\draw [line width=2pt] (0.5,11)-- (-0.5,12);
\draw [color=zzttgg, line width=1pt] (-0.5,12)-- (-1.2,12.7);
\draw (-1.35,12.4)		node[anchor=south] {\small\color{green}b};
\draw (-0.5,12)		node[anchor=east] {\small5};
\draw [color=zzttgg, line width=1pt] (-0.5,12)-- (0.2,12.7);
\draw (0.35,12.4)		node[anchor=south] {\small\color{green}c};
\draw (-0.5,12)		node[anchor=south] {\small$\overline{1}$};
\draw [color=zzttgg, line width=1pt] (0.5,11)-- (1.2,11.7);
\draw (1.35,11.4)		node[anchor=south] {\small\color{green}d};
\draw (0.5,11)		node[anchor=south] {\small$\overline{2}$};
\draw [color=zzttgg, line width=1pt] (1.5,10)-- (2.2,10.7);
\draw (2.35,10.4)		node[anchor=south] {\small\color{green}e};
\draw (1.5,10)		node[anchor=south] {\small$\overline{2}$};
\draw [color=zzttgg, line width=1pt] (2.5,9)-- (3.2,9.7);
\draw (3.35,9.4)		node[anchor=south] {\small\color{green}f};
\draw (2.5,9)		node[anchor=south] {\small$\overline{2}$};
\draw [color=zzttgg, line width=1pt] (3.5,8)-- (4.2,8.7);
\draw (4.35,8.4)		node[anchor=south] {\small\color{green}g};
\draw (3.5,8)		node[anchor=south] {\small$\overline{2}$};
\draw (3.5,8)		node[anchor=west] {\small$\overline{1}$};
\end{scope}

\begin{scope}[xshift=6cm]

\draw (3.5,7)		node[anchor=north] {\small\color{green}\bf a};
\draw [line width=2pt] (3.5,7)-- (3.5,8);
\draw [line width=2pt] (3.5,8)-- (2.5,9);
\draw [line width=2pt] (2.5,9)-- (1.5,10);
\draw [color=zzttgg, line width=1pt] (1.5,10)-- (0.8,10.7);
\draw (0.65,10.4)		node[anchor=south] {\small\color{green}b};
\draw (1.5,10)		node[anchor=east] {\small3};
\draw [line width=2pt] (1.5,10)-- (2.5,11);
\draw [color=zzttgg, line width=1pt] (2.5,11)-- (1.8,11.7);
\draw (1.65,11.4)		node[anchor=south] {\small\color{green}c};
\draw (2.5,11)		node[anchor=east] {\small2};
\draw [color=zzttgg, line width=1pt] (2.5,11)-- (3.2,11.7);
\draw (3.35,11.4)		node[anchor=south] {\small\color{green}d};
\draw (2.5,11)		node[anchor=south] {\small$\overline{1}$};
\draw [color=zzttgg, line width=1pt] (2.5,9)-- (3.2,9.7);
\draw (3.35,9.4)		node[anchor=south] {\small\color{green}e};
\draw (2.5,9)		node[anchor=south] {\small$\overline{3}$};
\draw [line width=2pt] (3.5,8)-- (4.5,9);
\draw [color=zzttgg, line width=1pt] (4.5,9)-- (3.8,9.7);
\draw (3.65,9.4)		node[anchor=south] {\small\color{green}f};
\draw (4.5,9)		node[anchor=east] {\small3};
\draw [color=zzttgg, line width=1pt] (4.5,9)-- (5.2,9.7);
\draw (5.35,9.4)		node[anchor=south] {\small\color{green}g};
\draw (4.5,9)		node[anchor=south] {\small$\overline{1}$};
\draw (3.5,8)		node[anchor=west] {\small$\overline{2}$};

\end{scope}

\end{tikzpicture}\\
\begin{tikzpicture}[scale=.6]

\begin{scope}[xshift=0cm]
\polygon[scale=2,xshift=5cm,yshift=5cm]{7}{{1/3,1/4,1/5,1/6}}
\draw [color=zzttgg, line cap=round, line join=round, line width=1pt] (8.540,9.667)-- (8.515,8.815);
\draw [color=zzttgg, line cap=round, line join=round, line width=1pt] (8.540,9.667)-- (8.148,10.423);
\draw [line cap=round, line join=round, line width=2pt] (9.190,10.482)-- (8.540,9.667);
\draw [color=zzttgg, line cap=round, line join=round, line width=1pt] (9.190,10.482)-- (9.176,11.712);
\draw [line cap=round, line join=round, line width=2pt] (10.232,10.482)-- (9.190,10.482);
\draw [color=zzttgg, line cap=round, line join=round, line width=1pt] (10.232,10.482)-- (10.824,11.712);
\draw [line cap=round, line join=round, line width=2pt] (10.882,9.667)-- (10.232,10.482);
\draw [color=zzttgg, line cap=round, line join=round, line width=1pt] (10.882,9.667)-- (11.852,10.423);
\draw [line cap=round, line join=round, line width=2pt] (10.650,8.650)-- (10.882,9.667);
\draw [color=zzttgg, line cap=round, line join=round, line width=1pt] (10.650,8.650)-- (11.485,8.815);
\draw [line cap=round, line join=round, line width=2pt] (10.650,8.650)-- (10.000,8.000);
\end{scope}

\begin{scope}[xshift=9cm]

\polygon[scale=2,xshift=5cm,yshift=5cm]{7}{{2/4,1/4,1/5,5/7}}
\draw [color=zzttgg, line cap=round, line join=round, line width=1pt] (8.829,10.934)-- (8.148,10.423);
\draw [color=zzttgg, line cap=round, line join=round, line width=1pt] (8.829,10.934)-- (9.176,11.712);
\draw [color=zzttgg, line cap=round, line join=round, line width=1pt] (9.061,9.918)-- (8.515,8.815);
\draw [line cap=round, line join=round, line width=2pt] (9.061,9.918)-- (8.829,10.934);
\draw [line cap=round, line join=round, line width=2pt] (10.232,10.482)-- (9.061,9.918);
\draw [color=zzttgg, line cap=round, line join=round, line width=1pt] (10.232,10.482)-- (10.824,11.712);
\draw [color=zzttgg, line cap=round, line join=round, line width=1pt] (11.460,9.667)-- (11.852,10.423);
\draw [color=zzttgg, line cap=round, line join=round, line width=1pt] (11.460,9.667)-- (11.485,8.815);
\draw [line cap=round, line join=round, line width=2pt] (10.521,9.214)-- (10.232,10.482);
\draw [line cap=round, line join=round, line width=2pt] (10.521,9.214)-- (11.460,9.667);
\draw [line cap=round, line join=round, line width=2pt] (10.521,9.214)-- (10.000,8.000);

\end{scope}

\end{tikzpicture}\\

The reasoning above can be extended in order to obtain -- in a recursive manner -- descriptions of similarity types (of quiddity sequences of arbitrary length).

We will sketch how this is done. Consider  $\bf t \it = type(\bf a\rm )$ a similarity type given by a quiddity sequence $\bf a\it =(a_{\rm 0\it},a_{\rm 1\it},...,a_{n-\rm 1\it })$ of length $n=len(\bf a \rm)$. In order to recover all the quiddity sequences of this type $\bf t \rm$,  it is necessary to distinguish between three categories of quiddity sequences, and to consider their periods.

The \bf period \rm $p=per(\bf a\rm )$ of a quiddity sequence $\bf a\it =(a_{\rm 0\it},a_{\rm 1\it},...,a_{n-\rm 1\it })$ of length $n$ is, obviously, the minimal strictly positive number satisfying the condition $a_p=a_0$. Since $a_n=a_0$, it is clear that $n$ is a multiple of $p$; in fact, $n=p$, $n=2p$ or $n=3p$ (see \cite{Cuntz}).

The \bf category \rm of the  quiddity sequence $\bf a\rm $ depends on its period $p$. If $p$ is odd and  $(a_0,a_1,...,a_{p-1})$ is symmetric, i.e. $a_i=a_{p-i-1}$ for $i\in \{ 0,...,p-1\} $, then the quiddity sequence $\bf a\rm $ is called \bf symmetric \rm (see above). If $p$ is even and  $a_i=a_{p-i}$ for $i\in \{ 1,...,p-1\} $, we say that $\bf a\rm $ is \bf pseudo-symmetric\rm . All the other quiddity sequences will be \bf asymmetric\rm .

\vspace{6pt}
Now, let $n\ge 3$ and denote by $\mathcal{K}_n$ the family of similarity types of quiddity sequences of length $n$. Of course, any such similarity type $\bf t \rm$ is represented by (at least) a quiddity sequence $\bf a\it =(a_{\rm 0\it},a_{\rm 1\it},...,a_{n-\rm 1\it})$. Notation $\bf t \it = type(\bf a\rm )$. 

Thus, $\mathcal{K}_3=\{ type(1,1,1)\} $ (of period 1), $\mathcal{K}_4=\{ type(2,1,2,1)\} $ (of period 2) and $\mathcal{K}_5=\{ type(2,1,3,1,2)\} $, the last quiddity sequence being a symmetric one. In addition, $\mathcal{K}_6$ can be described as 

$\{ type(3,1,2,3,1,2), type(4,1,2,2,2,1), type(3,1,3,1,3,1)\} $,

\noindent the last two being pseudo-symmetric.

If $\bf a\rm $ is a quiddity sequence of length $n$ that represents the type $\bf t \rm $ and $p=per(\bf a\rm )$, then every quiddity sequence $\delta *\bf a\rm $, where $\delta\in\mathbb{D}_p$, represents the same similarity type, i.e. $type(\delta *\bf a\rm )=\bf t $\rm. Moreover, all the quiddity sequences that are similar to $\bf a\rm $ can be recovered as $\delta*\bf a\rm $ where $\delta \in \mathbb{D}_p$.

(In the notatios above, symmetric quiddity sequences $\bf a\rm$ are those satisfying the condition $\sigma*\bf a \it = \bf a \rm$, and pseudo-symmetric ones satisfy $\sigma*\bf a \it = \gamma*\bf a \rm$.)

\vspace{6pt}
Now, consider a given positive integer $n$ and a (perfect) partition $(i,j,k)$ of it. Once all the quiddity sequences of lengths $i+1$, $j+1$, $k+1$ have been recovered from the respective types of  $\mathcal{K}_{i+1}$, $\mathcal{K}_{j+1}$ and $\mathcal{K}_{k+1}$, they will be composed, according to the following ternary composition law.

This law is defined for three arbitrary quiddity sequences $\bf a\it =(a_{\rm 0\it},a_{\rm 1\it},...,a_{u})$, $\bf b\it=(b_{\rm 0\it},b_{\rm 1\it},...,b_{v})$ and $\bf c\it=(c_{\rm 0\it},c_{\rm 1\it},...,c_{w})$, as follows 

\vspace{3pt}
$\langle \bf a, b, c\it \rangle = \bf z\it$

\vspace{3pt}
\noindent where $\bf z\it =(z_{\rm 0\it},z_{\rm 1\it},...,z_{u+v+w-\rm 1})$ is the quiddity sequence of length $u+v+w$ whose components are: $z_0=a_0+c_w$, $z_i=a_i$ for $1\le i<u$, $z_{u}=a_u+b_0$, $z_{u+j}=b_{j}$ for $1\le j<v$, $z_{u+v}=b_v+c_0$, 
and $z_{u+v+k}=c_k$ for $1\le k<w$.
 
 $\mathcal{K}_{n}$ will eventually be obtained by assembling all these compositions for all the perfect tri-partitions of $n$. 

(However, the above construction is rather a theoretical one. It is difficult to imagine how the list  $\mathcal{K}_{16}$, containing 83898 types, is completed.)

\end{document}